Abramov E. G.[1]

# Unconstrained inverse quadratic programming problem

**Abstract**. The paper covers a formulation of the inverse quadratic programming problem in terms of unconstrained optimization where it is required to find the unknown parameters (the matrix of the quadratic form and the vector of the quasi-linear part of the quadratic form) provided that approximate estimates of the optimal solution of the direct problem and those of the target function to be minimized in the form of pairs of values lying in the corresponding neighborhoods are only known. The formulation of the inverse problem and its solution are based on the least squares method. In the explicit form the inverse problem solution has been derived in the form a system of linear equations. The parameters obtained can be used for reconstruction of the direct quadratic programming problem and determination of the optimal solution and the extreme value of the target function, which were not known formerly. It is possible this approach opens new ways in over applications, for example, in neurocomputing and quadric surfaces fitting. Simple numerical examples have been demonstrated. A scenario in the Octave/MATLAB programming language has been proposed for practical implementation of the method.

## 1. Introduction

Inverse optimization problems are known, i.e. such problems where it is required to find (or select) precisely the values of parameters converting a certain desirable permissible solution into the optimal one. Inverse problems of discrete optimization have already been known sufficiently well, e.g. a good survey can be found in C. Heuberger [1]. Inverse problems of continuous optimization have been poorly studied so far, though, for example, in Iyengar and Kang [2] their practical relevance has been shown by example of inverse conic programming problems. However, for the first time, continuous inverse optimization problems by example of linear programming were formulated and studied by J. Zhang and Liu [3, 4]. A considerable time afterwards, at first J. Zhang et al. [5, 6], then also other authors [7 - 9] formulated inverse quadratic programming problems. All these formulations of inverse quadratic programming problems have a common feature – a certain approximate estimate of parameters is known beforehand, though it does not make the permissible solution to be optimal. This paper proposes another formulation of the inverse quadratic programming problem implying that the parameters are unknown at all but the pairs of values lying in the neighborhoods of permissible solution and the target function corresponding thereto are approximately known.

Let us consider a direct quadratic programming problem:

[1] evgeniy.g.abramov@gmail.com
Saint Petersburg State University
7/9 Universitetskaya emb., St. Petersburg, 199034, Russia
Submitted on 04 Jan 2017(v1)

$$\min_{\mathbf{x}} f(\mathbf{x}) = \frac{1}{2}\mathbf{x}^T\mathbf{G}\mathbf{x} + \mathbf{c}^T\mathbf{x} \qquad (1.1)$$
$$\text{subject to } \mathbf{A}\mathbf{x} \leq \mathbf{b}$$

In most general terms, J. Zhang et al. [6] offer the following formulation of the inverse problem.

Suppose that a solution $\mathbf{x}_0$ is known beforehand with a certain approximate estimate of parameters ($\mathbf{G}_0$, $\mathbf{c}_0$, $\mathbf{b}_0$), then the search for the precise values ($\mathbf{G}$, $\mathbf{c}$, $\mathbf{b}$) meeting (1.1) can be formulated as the following inverse problem:

$$\min_{\mathbf{x}} \|(\mathbf{G}, \mathbf{c}, \mathbf{b}) - (\mathbf{G}_0, \mathbf{c}_0, \mathbf{b}_0)\|^2 \qquad (1.2)$$
$$\text{where } \mathbf{x_0} \in \text{SOL}(QP(G, c, A, b))$$

The solutions of problems of type (1.2) require engaging the Karush-Kuhn-Tucker conditions.

This paper proposes another problem being inverse with respect to (1.1) reducing to unconstrained optimization.

That is to say:

**Problem formulation**

*Suppose that in problem (1.1) neither parameters (**G**, **c**), no precise solution $\mathbf{x}_0$, no minimal value of the target function $f(\mathbf{x}_0)$ corresponding thereto are known. However at that the pairs of values $(\mathbf{x}_i, y_i)$ lying accordingly in the neighborhoods $\mathbf{x}_0$ and $f(\mathbf{x}_0)$ are known, so that if the (**G**, **c**) were known, then the following equalities would be approximately (or even precisely) valid:*

$$y_i \approx \frac{1}{2}\mathbf{x}_i^T\mathbf{G}\mathbf{x}_i + \mathbf{c}^T\mathbf{x}_i \qquad (1.3)$$

*Then let us formulate the problem being inverse with respect to (1.1) as a search for the unknown values (**G**, **c**), basing on the approximate (or precise) estimates of $\mathbf{x}_0$ and $f(\mathbf{x}_0)$ in the form of the pairs of values $(\mathbf{x}_i, y_i)$, as follows:*

$$\min_{\mathbf{G},\mathbf{c}} \Phi(\mathbf{G}, \mathbf{c}) = \frac{1}{2}\sum_{i=1}^{N}\left(\frac{1}{2}\mathbf{x}_i^T\mathbf{G}\mathbf{x}_i + \mathbf{c}^T\mathbf{x}_i - y_i\right)^2 \qquad (1.4)$$

## 2. Main result: solution of the inverse problem

Let is introduce a constant multiplier to (1.4) so that as if both parts of equalities (1.3) were multiplied by 2:

$$\min_{\mathbf{G},\mathbf{c}} \Phi(\mathbf{G}, \mathbf{c}) = \frac{1}{2}\sum_{i=1}^{N}\left(\mathbf{x}_i^T\mathbf{G}\mathbf{x}_i + 2\mathbf{c}^T\mathbf{x}_i - 2y_i\right)^2 \qquad (2.1)$$

Then let us substitute variables as follows:

$$\hat{\mathbf{x}}_i = [1, \mathbf{x}_i^T]^T, \text{where } \hat{x}_{0i} = 1, \hat{x}_{li} = x_{li}, \forall l = 1, \cdots, m, \forall i = 1, \cdots, N \qquad (2.2)$$

$$\mathbf{W} = \begin{pmatrix} w_{00} & \mathbf{c}^T \\ \mathbf{c} & \mathbf{G} \end{pmatrix}, \text{where } w_{p0} = w_{0p} = c_p, w_{pr} = g_{pr}, \forall p = 1, \cdots, m, \forall r = 1, \cdots, m \qquad (2.3)$$

After that we derive formulation (1.4) in the simplified notation:

$$\min_{\mathbf{W}} Q(\mathbf{W}) = \frac{1}{2}\sum_{i=1}^{N}(\hat{\mathbf{x}}_i^T\mathbf{W}\hat{\mathbf{x}}_i - 2y_i)^2 \qquad (2.4)$$

Let us differentiate the function $Q(\mathbf{W})$ with respect to $\mathbf{W}$ and equalize the result to zero:

$$\frac{\partial Q}{\partial w_{pr}} = \sum_{i=1}^{N}(\hat{\mathbf{x}}_i^T \mathbf{W} \hat{\mathbf{x}}_i - 2y_i)\hat{x}_{pi}\hat{x}_{ri} = 0, \forall p = 0, \cdots, m, \forall r = 0, \cdots, m \quad (2.5)$$

Let us rewrite (2.5) in the scalar form having simultaneously performed regrouping of the multipliers:

$$\sum_{i=1}^{N}\sum_{l=0}^{m}\sum_{t=0}^{Z} \hat{x}_{li} w_{lt} \hat{x}_{ti} \hat{x}_{pi} \hat{x}_{ri} = 2 \sum_{i=1}^{N} y_i \hat{x}_{pi} \hat{x}_{ri}, \forall p = 0, \cdots, m, \forall r = 0, \cdots, m \quad (2.6)$$

In essence, $(m+1)^2$ equations of type (2.6) form a square non-homogenous system of linear (with respect to unknown variables $w_{lt}$) equations. This system can be rewritten in the matrix form:

$$\begin{pmatrix} \left(\text{vec}\left(\widehat{\mathbf{X}} \mathbf{D}_{\hat{x}0} (\widehat{\mathbf{X}} \mathbf{D}_{\hat{x}0})^T\right)\right)^T \\ \cdots \\ \left(\text{vec}\left(\widehat{\mathbf{X}} \mathbf{D}_{\hat{x}p} (\widehat{\mathbf{X}} \mathbf{D}_{\hat{x}r})^T\right)\right)^T \\ \cdots \\ \left(\text{vec}\left(\widehat{\mathbf{X}} \mathbf{D}_{\hat{x}m} (\widehat{\mathbf{X}} \mathbf{D}_{\hat{x}m})^T\right)\right)^T \end{pmatrix} \text{vec}\mathbf{W} = \text{vec}(2\widehat{\mathbf{X}} \mathbf{D}_y \widehat{\mathbf{X}}^T), \quad (2.7)$$

where $\widehat{\mathbf{X}} = [\hat{x}_{pi}]_{(m+1)\times N}$, $\mathbf{D}_{\hat{x}p} = \text{diag}(\hat{x}_{p1}, \hat{x}_{p2}, \cdots, \hat{x}_{pi}, \cdots, \hat{x}_{pN}), \forall p = 0, \cdots, m$; $\mathbf{D}_{\hat{x}r} = \text{diag}(\hat{x}_{r1}, \hat{x}_{r2}, \cdots, \hat{x}_{ri}, \cdots, \hat{x}_{rN}), \forall r = 0, \cdots, m$; $\mathbf{D}_y = \text{diag}(y_1, y_2, \cdots, y_i, \cdots, y_N)$, **vec** defines vectorization of a matrix.

Further determination of the solution is in principle an easy job. However difficulties of computational character may arise.

## 3. Discussion

As can be seen from (2.7) the number of equations in the system depends on the dimension of the problem $m$ rather than on the number of the pairs of values $(\mathbf{x}_i, y_i)$. Additional investigations of this system is required for the requirements to be made for the input data ($\mathbf{X}$ and $\mathbf{y}$) in order to ensure the system compatibility conditions according to the *Kronecker-Capelli theorem* as well as its well-conditionally. In particular, it is important to answer the following **question**: *how the condition number of system (2.7) varies when the number of pairs of values $(\mathbf{x}_i, y_i)$ changes in respect to $(m+1)^2$*. Generally, in linear algebra or regression handbooks (for example, [10] and [11]) the matters of ill-conditioned linear systems and methods for dealing with it such as *Tikhonov regularization*, *pseudoinversion*, *preconditioning* are considered in separate chapters. At the current stage of investigation, system (2.7) is solved using the *Moor-Penrose pseudoinverse*. **Appendix A** contains a scenario in the Octave/MATLAB language using this pseudoinverse for determination of matrix $\mathbf{W}$.

It should be noted that the solution of problem (1.4) allows to completely reconstruct the direct problem (1.1) knowing the additional constraints ($\mathbf{A}$, $\mathbf{b}$), i.e. determine not only the parameters ($\mathbf{G}$, $\mathbf{c}$) but also the initially unknown solution $\mathbf{x}_0$, and the minimum value of the target function $f(\mathbf{x}_0)$ corresponding thereto.

Really, unconstrained inverse quadratic programming problem (UIQPP) is only most evident application of main result. It is possible this approach opens new ways in over applications, for example, in neurocomputing [12] and quadric surfaces fitting problems [13 - 15]. Let us consider simple computational examples.

## 4. Simple examples of computations

Let us consider two low-dimensional direct convex quadratic programming problems (1.1) with the known set of parameters:

$$\mathbf{G}_1 = \begin{pmatrix} 2 & 1 \\ 1 & 2 \end{pmatrix}, \mathbf{c}_1 = \begin{pmatrix} 1 \\ 2 \end{pmatrix}, \mathbf{A}_1 = \begin{pmatrix} 1 & 2 \\ 1 & 3 \end{pmatrix}, \mathbf{b}_1 = \begin{pmatrix} -3 \\ -4 \end{pmatrix} \quad (4.1)$$

$$\mathbf{G}_2 = \begin{pmatrix} 4 & 1 & 2 \\ 1 & 4 & 3 \\ 2 & 3 & 4 \end{pmatrix}, \mathbf{c}_2 = \begin{pmatrix} 1 \\ 2 \\ 3 \end{pmatrix}, \mathbf{A}_2 = \begin{pmatrix} 1 & 2 & 3 \\ 1 & 3 & 4 \\ 1 & 4 & 5 \end{pmatrix}, \mathbf{b}_2 = \begin{pmatrix} -3 \\ -4 \\ -5 \end{pmatrix} \quad (4.2)$$

Then the optimal solutions of problem (1.1) with the corresponding values of target functions will be as follows:

$$\mathbf{x}_0^{(1)} = \begin{pmatrix} 0 \\ -1.5 \end{pmatrix}, f^{(1)}\left(\mathbf{x}_0^{(1)}\right) = -0.75 \quad (4.3)$$

$$\mathbf{x}_0^{(2)} = \begin{pmatrix} 0.25 \\ 0.2499 \\ -1.25 \end{pmatrix}, f^{(2)}\left(\mathbf{x}_0^{(2)}\right) = -1.125, \quad (4.4)$$

### 4.1. A case with strict equalities

Let us suppose that in the inverse problem (1.4) the values in the neighborhoods of $\mathbf{x}_0$ and $f(\mathbf{x}_0)$ are known so well that the approximate equalities (1.3) would inverse into the strict ones if only parameters $(\mathbf{G}, \mathbf{c})$ were known:

$$y_i = \frac{1}{2}\mathbf{x}_i^T \mathbf{G} \mathbf{x}_i + \mathbf{c}^T \mathbf{x}_i \quad (4.1.1)$$

As an example, we will be guided by the following values of pairs $(\mathbf{x}_i, y_i)$ from the neighborhoods of values (4.3) and (4.4):

$$\mathbf{X}_1 = \begin{pmatrix} 0 & 0.1 & 0.1 & 0.1 & 0 & 0.2 & 0.2 & -0.1 & -0.1 & 0 & 0.2 & 0.4 & -0.1 & -0.1 & -0.2 & -0.2 & -0.2 & -0.2 & 0 & 0.3 \\ -1.8 & -1.8 & -1.7 & -1.6 & -1.7 & -1.6 & -1.7 & -1.5 & -1.6 & -1.6 & -1.8 & -1.8 & -1.7 & -1.8 & -1.8 & -1.7 & -1.6 & -1.5 & -1.5 & -1.7 \end{pmatrix} \quad (4.1.2)$$

$$\mathbf{y}_1^T = (-0.36 - 0.43 - 0.57 - 0.69 - 0.51 - 0.72 - 0.61 - 0.69 - 0.57 - 0.64 - 0.48 - 0.52 - 0.43 - 0.27 - 0.16 - 0.33 - 0.48 - 0.61 - 0.75 - 0.63)$$

$$\mathbf{X}_2 = \begin{pmatrix} 0.2 & 0.1 & 0.2 & 0.3 & 0.1 & 0 & 0.3 & 0 & 0.3 & 0 & 0.2 & 0 & 0.1 & 0 & 0 & 0.2 & 0 & 0 & 0.2 & 0 \\ 0.2 & 0.1 & 0.2 & 0.3 & 0 & 0.1 & 0 & 0.3 & 0 & 0 & 0 & 0.1 & 0 & 0.1 & 0 & 0 & 0.2 & 0 & 0 & 0.2 \\ -1.3 & -1.4 & -1.4 & -1.4 & -1.3 & -1.3 & -1.3 & -1.3 & -1.2 & -1.1 & -1.1 & -1.1 & -1.2 & -1.2 & -1.2 & -1.2 & -1.2 & -1.3 & -1.3 & -1.3 \end{pmatrix} \quad (4.1.3)$$

$$\mathbf{y}_2^T = (-1.02 - 0.63 - 0.88 - 1.03 - 0.66 - 0.69 - 0.82 - 0.91 - 0.96 - 0.88 - 1.04 - 0.99 - 0.84 - 0.86 - 0.72 - 0.92 - 0.96 - 0.52 - 0.76 - 0.82)$$

Solving system (2.7) with respect to $\mathbf{W}$, we derive the following values:

$$\mathbf{W}_1 = \begin{pmatrix} 0.0000 & 1.0000 & 2.0000 \\ 1.0000 & 2.0000 & 1.0000 \\ 2.0000 & 1.0000 & 2.0000 \end{pmatrix} \quad (4.1.4)$$

$$\mathbf{W}_2 = \begin{pmatrix} 0.0000 & 1.0000 & 2.0000 & 3.0000 \\ 1.0000 & 4.0000 & 1.0000 & 2.0000 \\ 2.0000 & 1.0000 & 4.0000 & 3.0000 \\ 3.0000 & 2.0000 & 3.0000 & 4.0000 \end{pmatrix} \quad (4.1.5)$$

which fully conforms to the initial parameters $(\mathbf{G}, \mathbf{c})$ in the direct problem (1.1).

### 4.2. A case with approximate equalities

Now let us suppose that in the inverse problem (1.4) the values in the neighborhoods of $\mathbf{x}_0$ and $f(\mathbf{x}_0)$ are known with a certain error so that the approximate equalities (1.3) are valid. As an example, we will be guided by distorted values of pars $(\mathbf{x}_i, y_i)$ (4.1.2) and (4.1.3) from the neighborhoods of values (4.3) and (4.4):

$$X_1 = \begin{pmatrix} 0.01 & 0.11 & 0.11 & 0.11 & 0 & 0.19 & 0.21 & -0.09 & -0.1 & 0.01 & 0.2 & 0.41 & -0.09 & -0.09 & -0.2 & -0.2 & -0.19 & -0.19 & -0.01 & 0.3 \\ -1.8 & -1.79 & -1.69 & -1.6 & -1.7 & -1.6 & -1.7 & -1.5 & -1.6 & -1.6 & -1.79 & -1.79 & -1.7 & -1.79 & -1.79 & -1.7 & -1.59 & -1.5 & -1.5 & -1.7 \end{pmatrix} \quad (4.2.1)$$

$$y_1^T = (-0.36\ -0.42\ -0.57\ -0.68\ -0.51\ -0.71\ -0.61\ -0.68\ -0.56\ -0.63\ -0.48\ -0.52\ -0.43\ -0.26\ -0.16\ -0.32\ -0.47\ -0.6\ -0.75\ -0.63)$$

$$X_2 = \begin{pmatrix} 0.21 & 0.11 & 0.2 & 0.3 & 0.1 & 0 & 0.3 & 0.01 & 0.3 & 0 & 0.2 & 0.01 & 0.11 & 0.01 & 0.01 & 0.2 & 0 & 0 & 0.21 & 0 \\ 0.2 & 0.11 & 0.21 & 0.3 & 0 & 0.11 & 0 & 0.31 & 0.01 & 0 & 0 & 0.1 & 0.01 & 0.11 & 0 & 0 & 0.2 & 0.01 & 0.01 & 0.21 \\ -1.3 & -1.39 & -1.4 & -1.39 & -1.3 & -1.29 & -1.29 & -1.29 & -1.2 & -1.1 & -1.1 & -1.1 & -1.19 & -1.19 & -1.2 & -1.19 & -1.19 & -1.29 & -1.3 & -1.29 \end{pmatrix} \quad (4.2.2)$$

$$y_2^T = (-1.02\ -0.63\ -0.87\ -1.02\ -0.65\ -0.69\ -0.81\ -0.9\ -0.96\ -0.88\ -1.04\ -0.98\ -0.84\ -0.85\ -0.71\ -0.92\ -0.96\ -0.51\ -0.75\ -0.81)$$

Solving system (2.7) with respect to **W**, we derive the following values:

$$W_1 = \begin{pmatrix} 0.1843 & 0.6311 & 2.1355 \\ 0.6311 & 1.9466 & 0.7979 \\ 2.1355 & 0.7979 & 2.1058 \end{pmatrix} \quad (4.2.3)$$

$$W_2 = \begin{pmatrix} -0.6346 & 0.9177 & 1.8917 & 2.5545 \\ 0.9177 & 4.3753 & 0.7568 & 2.0166 \\ 1.8917 & 0.7568 & 3.7942 & 2.8972 \\ 2.5545 & 2.0166 & 2.8972 & 3.7385 \end{pmatrix} \quad (4.2.4)$$

which approximately conform (in the given example with the accuracy in rounding to nearest integer) to the initial parameters(**G**, **c**) in the direct problem (1.1).

Now knowing the approximate parameters (**G**, **c**) let us try to solve already the direct quadratic programming problem with parameters (4.2.3) and (4.2.4) instead of the values that were previously used in (4.3) and (4.4) under the same additional constraints (**A**, **b**). So we hope to give an approximate estimate of $x_0$ and $f(x_0)$ basing only on the values $(x_i, y_i)$ from their supposed neighborhoods. After solving we obtain:

$$x_*^{(1)} = \begin{pmatrix} 0.0323 \\ -1.5162 \end{pmatrix}, f_*^{(1)}\left(x_*^{(1)}\right) = -0.8351 \quad (4.2.5)$$

$$x_*^{(2)} = \begin{pmatrix} 0.2573 \\ 0.2524 \\ -1.2540 \end{pmatrix}, f_*^{(2)}\left(x_*^{(2)}\right) = -0.8031, \quad (4.2.6)$$

We can compare the values derived to (4.3) and (4.4) determined in solving the initial direct problem.

## 5. Conclusion

So, having considered the typical formulation of the inverse quadratic programming problem implying the availability of a certain approximate estimate of the parameters, we have proceeded to another formulation where not only the parameters but also the precise estimates of the permissible optimal solution and corresponding values of the target function are unknown. For the inverse problem formulated, the solution has been determined in the form of a linear system, though requiring further studies. Problem solution examples are shown. A reader can apply the scenario in the Octave/MATLAB language published in **Appendix A** to some data.

## 6. Appendix A

This Application publishes a scenario in the Octave/MATLAB language using the *Moor-Penrose pseudoinverse* for determination of matrix **W**, according to the following input data: matrix **X** and vector **y** corresponding to the pairs of values $(x_i, y_i)$. As is accepted in linear algebra, the vectors are written in columns.

```
function [W] = UIQPP (X, y)
[m, N] = size (X);
capX = [(ones(1,N));X];
```

```
        LeftMatrix = zeros((m+1)*(m+1));
        for p = 1:(m+1)
            for r = 1:(m+1)

                MatrixForTVecInLeftMatrix = 
(capX*diag(capX(p,:)))*(capX*diag(capX(r,:)))';
                TVecInLeftMatrix = MatrixForTVecInLeftMatrix(:);
                LeftMatrix(((p-1)*(m+1))+r,:) = TVecInLeftMatrix;
            end
        end
        MatrixForRightVector = 2*capX*(diag(y))*capX';
        RightVector = MatrixForRightVector(:);
        VecW = pinv(LeftMatrix)*RightVector;
        W = vec2mat(VecW, m+1);
```

end